\newtheorem{thm}{Theorem}
\newtheorem{lem}{Lemma}
\newtheorem{df}{Definition}

\newtheorem{pred}{Proposition}
\newtheorem{utv}{Claim}

\newcommand{\dvo}{{\it Proof}}

\documentclass[12pt]{article}
\usepackage{amsfonts}
\title{Jacobian Conjecture and Nilpotent Mappings}\author{Vik.S. Kulikov
\thanks{Partly supported by RFFI
 (No. 96-01-00614) and INTAS (No. 96-0713).}}
\date{        }

\begin{document}
\maketitle
\begin{abstract}
We prove the equivalence of the Jacobian Conjecture (JC(n)) 
and the Conjecture on the cardinality of the set of fixed points of a polynomial nilpotent mapping (JN(n)) and prove 
a series of assertions confirming JN(n).
\end{abstract}

\section*{Introduction }

 Let $F: \mathbb C^n \to \mathbb C^n$ be a morphism of complex 
affine spaces. If we choose coordinates in the image and preimage, then $F$ is defined by $n$ polynomials 
 $Y_i=F_i(X_1,...,X_n)\in \mathbb C[X]=
\mathbb C[X_1,...,X_n]$ in $n$ variables. Put
$$F'(X)= \det JF(X),$$
where
$$ JF= \left( \frac{\partial F_i}{\partial X_j}
\right) _{1\leq i,j\leq n} $$
is the Jacobi martix of $F$.

It is well-known that $F$ is locally invertible in a neighborhood of some point $x\in \mathbb C^n$ (as an analytic 
mapping) if and only if $F'(x)\neq 0$. Therefore, "$F'(x)\neq 0$ for all $x\in \mathbb C^n$"
is a necessary condition for the invertibility of $F$.

The famous {\it Jacobian Conjecture} ($JC(n)$) claims: \newline
{\bf Conjecture} $JC(n)$. {\it Let $F: \mathbb C^n \to \mathbb C^n$ be a morphism such that $F'(x)\neq 0$ for all $x\in \mathbb C^n$ (or, equivalently, $F'(X)\in \mathbb C^{*}$), then $F$ is an isomorphism, i.e. $F$ possesses an inverse mapping which is  also given by polynomials.}

Without loss of generality we can assume that $F'(X)\equiv 1$ if $F'(X)\in \mathbb C^{*}$. The following reformulations of 
conjecture $JC(n)$ are well-known
\begin{thm} (cf. \cite{Bass}).
Let $F: \mathbb C^n \to \mathbb C^n$ be a morphism with $F'(X)\in \mathbb C^{*}$. The Jacobian conjecture $JC(n)$ is equivalent to either of the following assertions:

(i)  $F$ is injective;

(ii) The degree $\deg F$ of $F$ is equal to 1
 (i.e. $\mathbb C(X_1,...,X_n))=\mathbb C(F_1(X),...,F_n(X)) $);

(iii) The ring $\mathbb C[X_1,...,X_n]$ is a finitely generated
$\mathbb C[F_1(X),...,F_n(X)]$-module;

(iv) The ring $\mathbb C[X_1,...,X_n]$ is a projective
$\mathbb C[F_1(X),...,F_n(X)]$-module;

(v) The ring $\mathbb C[F_1(X),...,F_n(X)]$  is integrally closed in $\mathbb C[X_1,...,X_n]$;

(vi) The field $C(X_1,...,X_n)$ is a Galois extension over $\mathbb C(F_1(X),...,F_n(X))$;

(vii)  $F$ is proper.

\end{thm}

The Jacobian Conjecture have been formulated by O.H.Keller in 1939
\cite{Kel} in the case $n=2$ for the polynomials with integer coefficients.
For $n=1$ conjecture $JC(n)$ is trivial since a polynomial 
$F(X_1)\in \mathbb C[X_1]$ with $F'(X_1)\in \mathbb C^{*}$ is of 
degree one and, therefore, defines an isomorphism $F: \mathbb C^1 \to \mathbb C^1$.
The Jacobian Conjecture $JC(n)$ remains an open problem even in the case $n=2$ in spite of attempts of many mathematicians to prove it or to find a counterexample (a survey of the results on the Jacobian Conjecture one can find in  \cite{Bass}, \cite {Ess1}, \cite {Ess2}). These attempts allow to obtain many partial results connected with the Jacobian Conjecture and even several false proofs were published (for 
the critique of these false proofs we  refer to  \cite{Vit} and \cite{Bass} ).

The aim of this short note is to attract attention of the reader to 
a connection of $JC(n)$ and one of conjectures on the nilpotent mappings.

\begin{df}
A morphism $N=(N_1,...,N_n): \mathbb C^n \to \mathbb C^n$
is called nilpotent if
$JN(X)$ is a nilpotent matrix, i.e.
$$ (JN(X))^n= \left( \frac{\partial N_i}{\partial X_i}
\right)^n\equiv 0 $$
in the ring of matrices ${\cal{M}}(n, \, \mathbb C[X])$ of $n$-th order with coefficients in $\mathbb C[X]$.
\end{df}

We consider nilpotent mappings as endomorphisms of $\mathbb C^n$. A point $x\in \mathbb C^n$ is called a {\it fixed point}
of a nilpotent mapping $N$ if $N(x)=x$. It is well-known that if $N$ is linear and nilpotent, then $N^n(\mathbb C^n)$
consists of a unique point. In particular, a linear nilpotent mapping $N$
possesses only one fixed point.  \newline
{\bf Conjecture} $JN(n)$. {\it A nilpotent mapping $N: \mathbb C^n \to \mathbb C^n$ possesses at most one fixed point.}

The following theorem describes a connection between the 
conjectures  $JC(n)$ and $JN(n)$.
\begin{thm}\label{te}
The conjectures $JC(n)$ and $JN(n)$ are equivalent, i.e. $JC(n)$ is true for all $n$  if and only if $JN(n)$ is true for all $n$.
\end{thm}

In section 1, we recall some well-known results relating to the Jacobian Conjecture. In section 2, we prove Theorem \ref{te} and a series of assertions confirming  $JN(n)$.

This paper was written during my stay at the Max-Planck-Institut f\"{u}r Mathematik in Bonn. It is a pleasure to thank the Institut for its hospitality and financial support.

\section{Stable equivalence of mappings and the Jacobian Conjecture}
{\bf 1.1.} Let $F: \mathbb C^n \to \mathbb C^n$ and $G: \mathbb C^n \to \mathbb C^n$ be morphisms given by polynomials $F_1(X),...,F_n(X)$ and $G_1(X),...,G_n(X)$ respectively. We shall say that $F$ and $G$ are {\it equivalent} if there exist automorphisms $H, \, R: \mathbb C^n \to \mathbb C^n$ such that the following diagram is commutative
\\

\begin{picture}(0,0)(-45,3)
\put(136,4){$\mathbb C^n$}
\put(141,0){\vector(0,-1){20}}
\put(128,-12){$G$}
\put(136,-33){$\mathbb C^n$}
\put(155,-30){\vector(1,0){45}}
\put(170,-40){$R$}
\put(206,-33){$\mathbb C^n.$}
\put(205,4){$\mathbb C^n$}
\put(155,7){\vector(1,0){45}}
\put(170,10){$H$}
\put(209,0){\vector(0,-1){20}}
\put(214,-12){$F$}
\end{picture}
\\
\\
\\
\\
\\
In other words, $F$ and $G$ are equivalent if
$$R_i(G_1(X),...,G_n(X))=F_i(H_1(X),...,H_n(X))$$
for all $i$, where $H_1(X),...,H_n(X)$ and $R_1(X),...,R_n(X)$ are polynomials in
 $\mathbb C[X_1,...,X_n]$ defining the automorphisms $H$ and $R$.
\begin{df} Morphisms $F: \mathbb C^n \to \mathbb C^n$ and
$G: \mathbb C^m \to \mathbb C^m$, $n\leq m$, are called stably equivalent if $F\times Id : \mathbb C^m=\mathbb C^n\times \mathbb C^{m-n} \to
\mathbb C^n\times \mathbb C^{m-n}= \mathbb C^m$
and $G$ are equivalent.
\end{df}

The following assertion is a consequence of the definition.
\begin{utv}
If $F$ and $G$ are stably equivalent, then
$$\deg F= \deg G.$$
In particular, if $F'(X)\in \mathbb C^{*}$ for $F$ and
$F$ and $G$ are stably equivalent, then $F$ and $G$ are invertible (or not invertible) simultaneously.
\end{utv}
{\bf 1.2.} To avoid confusion with the degree of a morphism, we shall call   $p$-{\it degree} or
{\it polynomial degree} of a morphism $F=(F_1,...,F_n): \mathbb C^n \to
\mathbb C^n$ the maximum of the degrees of $F_i(X)$ defining $F$,
$$\deg _p F=\max \deg F_i(X).$$
Of course, $p$-degree depends on a choice of affine 
coordinates in the image and preimage.

A usual construction to simplify the form of polynomials defining a morphism is the following: one replaces a mapping by a stably equivalent one and performs linear and triangular automorphisms of the image and the preimage.
\begin{df}
An automorphism $H:\mathbb C^n \to \mathbb C^n$ is called triangular if, in some coordinate system, it can be given by polynomials $H_i(X_1,...,X_n), \, \, i=1,...,n,$ of the form
$$H_i(X_1,...,X_n)=X_i + H_{i}(X_1,...,X_{i-1}),$$
where $H_{i}(X_1,...,X_{i-1})\in \mathbb C[X_1,...,X_{i-1}]$ is a polynomial in variables $X_1,...,X_{i-1}$.
\end{df}

The following theorem is well-known.
\begin{thm} \label{red} (\cite{Bass}, \cite{Yag}).
For any morphism $F:\mathbb C^n\to \mathbb C^n$ there exists a stably equivalent morphism $G:\mathbb C^m\to \mathbb C^m$ such that
$\deg _p G \leq 3$.
\end{thm}
{\bf 1.3.} Let $F:\mathbb C^n\to \mathbb C^n$, given by 
polynomials $F_i(X)$, map the origin to the origin, i.e. $F(0)=0$. A morphism
$F_t:\mathbb C^{n+1}\to \mathbb C^{n+1}$ given by the polynomials
$$
\begin{array}{lll}
\widetilde F_i(X_1,...,X_n,T) & =  & \displaystyle \sum_{j=1}^{\deg F_i}T^{j-1}
F_{(j)i}(X_1,...,X_n), \, \, \, \, \, \, \, i=1,...,n, \\
\widetilde F_{n+1}(X_1,...,X_n,T) & = & T,
\end{array}
$$
is called {\it a blow up} of $F$, where
$F_{(j)i}(X)$ is the homogeneous component of degree $j$ of the polynomial $F_i(X)$.
It is easy to check the following assertion.
\begin{utv} Let $F:\mathbb C^n\to \mathbb C^n$ be a morphism such that $F(0)=0$, then

(i) $F'_t(X_1,...,X_n,T)=F'(TX_1,...,TX_n)$;

(ii) $F$, with $F'(X)\in \mathbb C^{*}$, is invertible if and only if $F_t$ is invertible.
\end{utv}

Geometrically, to obtain $F_t$ from $F$, one needs to perform  $\sigma$-processes with centers at the origins of the image and the preimage of the morphism $F\times Id: \mathbb C^{n+1}=\mathbb C^n\times \mathbb C^1
\to \mathbb C^n\times \mathbb C^1=\mathbb C^{n+1}$ and to restrict the obtained morphism onto coordinates neighborhoods corresponding the last coordinate $T$.

In \cite{Druz}, theorem \ref{red} was made more precise.
\begin{thm} \label{J3}
Let $F:\mathbb C^n\to \mathbb C^n$ be a morphism such that $F(0)=0$ and $F'(0)\neq 0$. Then the blow up $F_t$ of $F$ is stably equivalent to $G:\mathbb C^m\to \mathbb C^m$ defined by polynomials
$G_i(X)$ of the form
\begin{equation}
G_i(X_1,...,X_m)=X_i+G_{(3)i}(X_1,...,X_m), \label{F_3}
\end{equation}
where $G_{(3)i}(X_1,...,X_m)$ are homogeneous forms of degree 3.
\end{thm}

Since the space of the homogeneous polynomials of degree $k$ ¯is generated by $k$-th powers of the linear polynomials,
the following theorem is a simple consequence of Theorem \ref{J3}.  
\begin{thm}
Let $F:\mathbb C^n\to \mathbb C^n$ be given by polynomials
$$F_i(X_1,...,X_n)=X_i+F_{(3)i}(X_1,...,X_n).$$
Then $F$ is stably equivalent to $G:\mathbb C^m\to \mathbb C^m$ 
given by polynomials $G_i(X)$ of the form
\begin{equation}
G_i(X_1,...,X_m)=X_i+(\sum a_{i,j}X_j)^3. \label{F_3l}
\end{equation}
\end{thm}
{\bf 1.4.} The Jacobian conjecture was proved in \cite{Wri} and \cite{Hub} for the morphisms of the form (\ref{F_3}) with $m\leq 4$ and in \cite{Hub} for the morphisms of the form (\ref{F_3l}) with $m\leq 7$. Besides, in \cite{Wan} the following theorem was proved 
\begin{thm} \label{JC2}
The conjecture $JC(n)$ is true for $F:\mathbb C^n\to \mathbb C^n$ of polynomial degree $\deg _p F \leq 2.$
\end{thm}
This theorem will be obtained in the next section as a consequence of more general result (see Proposition 1).

Note also that the conjecture $JC(n)$ was proved in \cite{Lan}
in the following two cases: when $F$ is given by polynomials
$F_i(X)\in X_i\mathbb C[X]$ and when $F_i(X)=X_i+\lambda _i M_i(X)$, where
$M_i(X)$ are monomials. \newline
{\bf 1.5.} In the case $F:\mathbb C^2\to \mathbb C^2$ the Jacobian conjecture holds for all $F$  of polynomial degree $\deg _p F\leq 100$ \cite{Moh}. 
In \cite{Mag} (see also \cite{Ess2}),
it was proved that if the degrees of polynomials defining $F$ with $F'(X)\in \mathbb C^{*}$ are coprime, then $F$ is invertible and moreover in this case one of polynomials must be linear. In particular, the Jacobian conjecture holds if one of the polynomials defining $F$ is of prime degree. In \cite{App}, this result was generalized to the case when the degree of one of the polynomials defining $F$ is a product of two different prime integers.

In addition, in \cite{Or}, it was shown that if $F:\mathbb C^2\to
\mathbb C^2$ with $F'(X)\equiv 1$ is not invertible, then $\deg F\geq 4$. \newline
{\bf 1.6.} A morphism $F:\mathbb C^{n}\to \mathbb C^{n}$ induces a polynomial mapping $\widetilde F:\mathbb R^{2n}\to \mathbb R^{2n}$ given by
$$\widetilde F=(Re\, F_1, Im\, F_1,...,Re\, F_n, Im\, F_n),$$
where $F_i$ are polynomials defining $F$. It is easy to see that
$$\det J\widetilde F=\mid \det JF\mid ^2.$$
Hence, $\det JF\in \mathbb C^{*}$ if and only if
$\det J\widetilde F\in \mathbb R^{*}$. \newline
{\bf Conjecture} $RJC(n)$. {\it A polynomial mapping
$F:\mathbb R^{n}\to \mathbb R^{n}$ with $\det J F\in \mathbb R^{*}$ is injective.}

It is easy to show that if $RJC(n)$ is true, then
$F:\mathbb R^{n}\to \mathbb R^{n}$ with $\det J F\in \mathbb R^{*}$ is invertible and the following proposition holds.
\begin{pred}
The conjecture $JC(n)$ is true for all $n$ if and only if the conjecture $RJC(n)$ is true for all $n$.
\end{pred}

Consider $F:\mathbb R^{n}\to \mathbb R^{n}$ with $\det J F\in \mathbb R^{*}$. Without loss of generality, we can assume that $F(0)=0$ and $JF(0) =E$, where $E$ is the identity matrix. Write $F$ in the form
\begin{equation}
F=Id +F_{(2)}+...+ F_{(m)},  \label{od}
\end{equation}
where $F_{(i)}$ are morphisms defined by homogeneous polynomials of degree $i$.
\begin{df}
A polynomial mapping $F:\mathbb R^{n}\to \mathbb R^{n}$ of the form (\ref{od}) is called positive (resp. negative) if all non-zero coefficients of the $F_{(i)}$ are positive (resp. negative).
\end{df}
In \cite{Yu}, the following theorem has been proved.
\begin{thm} \label{yu}
(i) The conjecture $RJC(n)$ holds for all $n$ and all negative polynomial mappings $F:\mathbb R^{n}\to \mathbb R^{n}$.

(ii) It suffices to prove $RJC(n)$ for all positive polynomial mappings $F:\mathbb R^{n}\to \mathbb R^{n}$.
\end{thm}

It follows from Theorem \ref{yu} (i) (passing to the blow up of $F$) that

{\it The conjecture $RJC(n)$ holds for all $n$ and for all positive polynomials mappings $F:\mathbb R^{n}\to \mathbb R^{n}$ of the form $F=Id +F_{(2)}+F_{(4)}+...
+ F_{(2m)}$, and also for all polynomial mappings  $F$ of the form $F=Id +F_{(2)}-F_{(3)}+...
+(-1)^{m} F_{(m)}$, where $F_{(i)}$ are morphisms defined by homogeneous polynomials of degree $i$ with non-negative coefficients.}

\section{Nilpotent Mappings}
{\bf 2.1.} It is well-known that a matrix $A\in {\cal M} (n,\mathbb C)$ is nilpotent if and only if $\det (tE-A)\equiv
t^n$, where $E$ is the identity matrix. As a consequence we obtain the following Lemma.
\begin{lem} \label{lem} A morphism $N:\mathbb C^n\to \mathbb C^n$ given by polynomials $N_i(X)$ is nilpotent if and only if $\widetilde F=(Id-T\cdot N)\times Id:
\mathbb C^n\times \mathbb C^1\to \mathbb C^n\times \mathbb C^1$, given by the polynomials
$$
\begin{array}{rll}
\widetilde N_i(X_1,...,X_n,T)& = & X_i-TN_i(X_1,...,X_n) \hspace{1cm} \mbox{for}\, \, \,  i=1,\, ...\, ,\, n, \\
\widetilde F_{n+1}(X_1,...,X_n,T) & = & T \, ,
\end{array}
$$ 
has $\widetilde F'(X,T)\equiv 1$.
\end{lem}

If we apply lemma \ref{lem} to the restriction of 
$\widetilde F=Id-T\cdot N$ to the hypersurface
$\mathbb C^n = \{ \, T=1\, \} $ we obtain that for the nilpotent mapping $N$ the corresponding morphism
$F=Id-N:\mathbb C^n\to \mathbb C^n$
has $F'(X)\equiv 1$. Besides, if there exist $d$
distinct fixed points of $N$, then for
$F=Id-N$ the preimage of the origin consists of $d$ distinct points, and hence degree of $F=Id-N$ is not less than $d$. In particular, if a nilpotent mapping
$N:\mathbb C^n\to \mathbb C^n$ possesses $d>1$ fixed points,
then $F=Id-N:\mathbb C^n\to \mathbb C^n$ is not invertible.  \newline
\vspace{0.2cm}
{\bf 2.2.} Theorem \ref{te} is a consequence of the following theorem.
\begin{thm} \label{??}
Any morphism $F:\mathbb C^n\to \mathbb C^n$ of degree $\deg F=d$ with $F'(X)\equiv 1$ is stably equivalent to $\overline F$ of the form
$$\overline F=Id-N,$$
where $N$ is a nilpotent morphism possessing $d$ fixed points.

 Conversely, if a nilpotent mapping  $N:\mathbb C^n\to \mathbb C^n$ possesses 
$d$ distinct fixed points, then $\overline F=Id-N$ is of degree not less than $d$.
\end{thm}
\dvo . Let $F:\mathbb C^n\to \mathbb C^n$ be a morphism of degree $\deg F=d$ with $F'(X)\equiv 1$. Performing a translation and linear transformations, if necessary, we can assume that $F(0)=0$, the number of the preimages of the origin $\# \{ F^{-1}(0) \} = \deg F$, and $F$
is given by polynomials
\begin{equation}
F_i(X_1,...,X_n)= X_i - F_{(\geq 2)i}(X_1,...,X_n),
\end{equation}
where $F_{(\geq 2)i}(X_1,...,X_n)$ are polynomials containing no monomials of degree zero and one.

Consider the blow up $F_t:\mathbb C^{n+1}\to \mathbb C^{n+1}$ of $F$.
The morphism $F_t$ is defined by polynomials
$$
\begin{array}{lllll}
Y_i & = & \widetilde F_i(X_1,...,X_n,T) & =  & X_i-\sum_{j=2}^{\deg F_i}T^{j-1}
F_{(j)i}(X_1,...,X_n), \, \, \, \, \, \, \, i=1,...,n, \\
\widetilde T & = & \widetilde F_{n+1}(X_1,...,X_n,T) & = & T.
\end{array}
$$
Let us show that $F_t$ is stably equivalent to
$\widetilde G=(Id-T\cdot G)\times Id:\mathbb C^m\times \mathbb C^1
\to \mathbb C^m\times \mathbb C^1$ given by
$\widetilde Y_i=\widetilde G_i(X_1,...,X_m,T)=X_i-TG_i(X_1,...,X_m)$ for $i=1,...,m,$
and $\widetilde G_{m+1}(X_1,...,X_n,T)=T$, and of degree $\deg \widetilde G =
\deg F$, and the number of the preimage of the point with coordinates $(0,...,0,1)$ is
$\# \{ \widetilde G^{-1}((0,...,0,1))\} =\deg F$. For this we must get rid of all monomials of the form $T^kX_1^{r_1}\cdots X_n^{r_n}$ with $k>1$ and then  apply Lemma
\ref{lem} to $G:\mathbb C^m\to \mathbb C^m$ defined by
$G_i(X_1,...,X_m)$, $i=1,...,m$. By Lemma 1, $G$ is nilpotent and the number of the fixed points of $G$ is equal to $\deg F$.

Removing the monomial $T^kX_1^{r_1}\cdots X_n^{r_n}$ with $k>1$ can be performed by adding new variables 
(passing to stably equivalent morphisms) and performing linear and triangular transformations. In fact, let for some $i$  a polynomial $\widetilde Y_i$ be of the form 
$$\widetilde Y_i=\widetilde X_i-\sum_{j=1}^{k_i}T^j F_{(j)i}(\widetilde X_1,...,\widetilde X_n)$$
(here $F_{(j)i}(\widetilde X_1,...,\widetilde X_n)$ is not necessary homogeneous). We add a new variable $\widetilde Y_{n+1}=\widetilde X_{n+1}$ and perform a triangular transformation in the preimage:
$$
\begin{array}{lllr}
X_i & = & \widetilde X_i, & i=1,...,n \\
X_{n+1} & = & \widetilde X_{n+1}+T F_{(k_i)i}(\widetilde X_1
,...,\widetilde X_n). &
\end{array}
$$
Next we perform a triangular transformation in the 
image:
$$
\begin{array}{llr}
Y_l= & \widetilde Y_l, & l=1,...,i-1,i+1,..., n+1 \\
Y_{i}= & \widetilde Y_{i}-T^{k_i-1}\widetilde Y_{n+1}. &
\end{array}
$$
It is easy to see that in the coordinates $X_1,...,X_{n+1},T$ and $Y_1,...,Y_{n+1},T$ the morphism $F_t\times Id$ is given by equations:
$$
\begin{array}{lllr}
Y_s & =  & X_s-\sum_{j=1}^{k_s}T^j F_{(j)s}(X_1,...,X_n),  & s=
1,...,i-1,i+1,..., n  \\
Y_i & = & X_i-\sum_{j=1}^{k_i-1}T^j F_{(j)l}(X_1,...,X_n)- T^{k_i-1}X_{n+1}, & \\
Y_{n+1} & = &  X_{n+1}-T F_{(k_i)i}(X_1,...,X_n). &
\end{array}
$$
Thus, we removed the summand of $\widetilde F_i(X_1,...,X_n,T)$ containing the variable $T$ is contained with the maximal multiplicity $k_i$. Obviously, after a finite number of similar steps we obtain a morphism stably equivalent to $F_t$ and satisfying the required conditions.
\newline  \vspace{0.5cm}
{\bf 2.3.} Let $F:\mathbb C^n\to \mathbb C^n$ be of the form
$F=Id-F_{(k)}$, i.e. $F$ is defined by
$$\mbox{\hspace{3cm}} F_i(X)=X_i -F_{(k)i}(X_1,...,X_n), \mbox{\hspace{3cm}}
i=1,...,n,$$
where $F_{(k)i}(X_1,...,X_n)$ are homogeneous forms of degree $k$. If $F'(X)\equiv 1$, then $F_{(k)}$ is nilpotent.
In fact, by Theorem 2, the Jacobian $F'_t(X,T)$ of the blow up $F_t$ is equal to 
$$F'_t(X,T)=F'(TX_1,...,TX_n)=\det ( E-T^{k-1}JF_{(k)}) \equiv 1.
$$
Therefore, the Jacobian of $\widehat F_t$ which is defined by polynomials
$$
\begin{array}{lllr}
Y_i & = & X_i-TF_{(k)i}(X_1,...,X_n), & i=1,...,n,  \\
T & = & T, &
\end{array}
$$
is also equal to 1. Hence, by Lemma \ref{lem}, $F_{(k)}$ is nilpotent.
\begin{pred}
The conjecture $JN(n)$ holds for nilpotent morphisms
$N_{(k)}:\mathbb C^n\to \mathbb C^n$ defined by homogeneous forms of degree $k$.
\end{pred}
\dvo . The origin, i.e. the point $o=(0,...,0)$,
 is a fixed point of $N=N_{(k)}$. Let $x=(x_1,...,x_n)\neq o$ be another fixed point of $N_{(k)}$. This contradicts the nilpotency of $N$, since the vector $\overline x=(x_1,...,x_n)$ is an eigenvector (corresponding to a non-zero eigenvalue (equals to $k$)) of the linear mapping $JN_{(k)}(X)$ at the point $X=x$. In fact, by Euler's formula
$$
\begin{array}{lll} JN_{(k)}(x)(\overline x) & = &
(\sum_{j=1}^{n}x_j\frac{\partial N_1}{\partial
X_j}(x),...,\sum_{j=1}^{n}x_j\frac{\partial N_n}{\partial
X_j}(x)) \\
& =  & (k N_1(x),...,k N_n(x)) = \\
& =  & (k x_1,...,k x_n) = k\cdot \overline x ,
\end{array}
$$
where $N_i(X)$  are polynomials defining $N=N_{(k)}$.

Theorem \ref{JC2} is an easy consequence of this theorem. For let $F$ be of polynomial degree $\deg _p F=2$. After a linear changes of variables, we can assume that $F$ is of the form $F=Id-F_{(2)}$ and the preimage  $F^{-1}(0)$ consists of $\deg F$ distinct points. Hence, the nilpotent morphism $F_{(2)}$ possesses $\deg F$ distinct fixed points. Therefore, $\deg F=1$, i.e. $F$ is invertible.

As a consequence of above considerations, Lemma \ref{lem}, Theorem \ref{JC2}, 
and the poof of Theorem \ref{??} we obtain the following theorem.
\begin{thm}
The conjecture $JN(n)$ holds for
$N:\mathbb C^n\to \mathbb C^n$ of polynomial degree $\deg _p N\leq 2$.
\end{thm}

Proposition 2 may be generalized as follows.
\begin{pred} \label{ff}
The conjecture $JN(n)$ holds for
$N:\mathbb C^n\to \mathbb C^n$ of the form $N=N_{(k_1)}+N_{(k_2)}$, where $N_{(k_i)}$ are morphisms defined by homogeneous forms of degree $k_i>0$.
\end{pred}
\dvo . The origin is a fixed point of $N$. Let $x=(x_1,...,x_n)\neq o$ be another fixed point of $N$. As in the proof of Proposition 2 this contradicts the nilpotency of $N$, since the vector $\overline x=(x_1,...,x_n)$ is an eigenvector (corresponding to a non-zero eigenvalue) of the linear mapping $JN(X)$ at the point $X=\lambda x=(\lambda x_1,...,\lambda x_n)$, where 
$$\lambda =\left( \frac{k_1}{k_2} \right) ^{\frac{1}{k_2-k_1}}.$$
In fact,
$$
\begin{array}{lll}
JN(\lambda x)(\overline x) & = &
(\sum_{j=1}^{n}x_j\frac{\partial N_{(k_1)1}}{\partial
X_j}(\lambda x),...,\sum_{j=1}^{n}x_j\frac{\partial N_{(k_1)n}}{\partial
X_j}(\lambda x))+ \\
 & + &(\sum_{j=1}^{n}x_j\frac{\partial N_{(k_2)1}}{\partial
X_j}(\lambda x),...,\sum_{j=1}^{n}x_j\frac{\partial N_{(k_2)n}}{\partial
X_j}(\lambda x)) = \\
& = & (k_1 \lambda ^{k_1-1}N_{(k_1)1}(x),...,k_1 \lambda ^{k_1-1}N_{(k_1)
n}(x))+ \\
& + & (k_2 \lambda ^{k_2-1}N_{(k_2)1}(x),...,k_2 \lambda ^{k_2-1}
N_{(k_2)n}(x))= \\
& =  &
k_1\left( \frac{k_1}{k_2} \right) ^{\frac{k_1-1}{k_2-k_1}}\cdot N_{(k_1)}(x)
+k_2\left( \frac{k_1}{k_2} \right) ^{\frac{k_2-1}{k_2-k_1}}\cdot N_{(k_2)}(x)
= 
\end{array}
$$

$$
\begin{array}{lll}
 & =  & \frac{k_1^{\frac{k_2-1}{k_2-k_1}}}{k_2^{\frac{k_1-1}{k_2-k_1}}}\cdot
N_{(k_1)}(x)+
\frac{k_1^{\frac{k_2-1}{k_2-k_1}}}{k_2^{\frac{k_1-1}{k_2-k_1}}}\cdot
N_{(k_2)}(x) \\
 & =  & \frac{k_1^{\frac{k_2-1}{k_2-k_1}}}{k_2^{\frac{k_1-1}{k_2-k_1}}}\cdot
(N_{(k_1)}(x)+ N_{(k_2)}(x)) \\
 & =  & \frac{k_1^{\frac{k_2-1}{k_2-k_1}}}{k_2^{\frac{k_1-1}{k_2-k_1}}}\cdot
\overline x .
\end{array}
$$
{\bf 2.4.} The following result is a consequence of Theorem \ref{red} and the proof of Theorem \ref{??} 

\begin{thm}
If each nilpotent morphism
$N:\mathbb C^n\to \mathbb C^n$ of polynomial degree $\deg _p N = 3$
has a unique fixed point, then $JN(n)$ holds.
\end{thm}

In the case $\deg _p N=3$ 
the morphism $N$ is of the form
$$N=N_{(1)}+N_{(2)}+N_{(3)},$$
where $N_{(i)}$ is a homogeneous form of degree $i$, and Proposition \ref{ff} can be reformulated as follows:

{\it The conjecture $JN(n)$ holds for
$N:\mathbb C^n\to \mathbb C^n$ of the form $N=N_{(1)}+N_{(3)}$ and} $N=N_{(2)}+N_{(3)}$.

It follows from \cite{Hub} that:

{\it The conjecture $JN(3)$ holds for
$N:\mathbb C^3\to \mathbb C^3$ of polynomial degree} $\deg _p N\leq 3$.
\newline
{\bf 2.5.} The following proposition is a consequence of Theorem \ref{yu}.
\begin{pred}
Let $F:\mathbb R^n \to \mathbb R^n$ be a polynomial mapping of the form
$F=Id+N$, where $N$ is a positive nilpotent morphism. Then $F$ is invertible.
\end{pred}
\dvo . Without loss of generality, we can assume that $N(0)=0$. Consider a morphism $\overline F:\mathbb R^{n+1} \to \mathbb R^{n+1}$ defined by
$\overline F_i=X_i-TN_i(X),\, \, i=1,...,n,$ and $\overline F_{n+1}=T$. By Theorem \ref{yu}, $\overline F$ is invertible. The restriction of $\overline F$ to the hyperplane
$T=-1$ coincides with $F$. Therefore, $F$ is invertible. 

In particular, it follows from this that:

{\it The conjecture $JN(n)$ holds for $N:\mathbb R^n \to \mathbb R^n$ if $N$ is either a positive or a negative nilpotent morphism. }

It follows from Theorem \ref{yu} (if we use a blow up) that:

{\it The conjecture $JN(n)$ holds for the nilpotent morphisms $N:\mathbb R^n \to \mathbb R^n$ of the form
$$N=\pm (N_{(1)}-N_{(2)}+...+(-1)^{m+1} N_{(m)}),$$
where $N_{(i)}$ are defined by positive homogeneous forms of degree $i$. }
\newline
{\bf 2.6.} Let $r= rk\, JN(X)$ be the rank of a nilpotent morphism $N:\mathbb C^n\to \mathbb C^n$. In \cite{Bass} it was shown that if $r=1$, then $F=Id-N$ is invertible:

{\it The conjecture $JN(n)$ holds for nilpotent morphisms
$N:\mathbb C^n\to \mathbb C^n$ of rank $r=1$. In particular, the conjecture $JN(2)$ is true.}

Steklov Mathematical Institute

{\rm victor$@$olya.ips.ras.ru }

{\rm kulikov$@$mpim-bonn.mpg.de}

\end{document}